\input amstex
\documentstyle{amsppt}
\magnification=\magstep1
\hcorrection{0in}
\vcorrection{0in}
\pagewidth{6.3truein}
\pageheight{9truein}
\topmatter
\title On Numerically Effective Log Canonical Divisors \endtitle
\author Shigetaka Fukuda
\leftheadtext{SHIGETAKA FUKUDA}
\rightheadtext\nofrills{\rm LOG CANONICAL DIVISORS}
\endauthor
\address Faculty of Education, Gifu Shotoku Gakuen University, Yanaizu-cho, Gifu 
501-6194, Japan
\endaddress
\email fukuda\@ha.shotoku.ac.jp
\endemail
\abstract
Let $(X,\Delta)$ be a 4-dimensional log variety which is proper over the field of 
complex numbers and with only divisorial log terminal singularities.
The log canonical divisor $K_X+\Delta$ is semi-ample, if it is nef (numerically 
effective) and the Iitaka dimension $\kappa(X,K_X+\Delta)$ is strictly positive.
For the proof, we use Fujino's abundance theorem for semi log canonical threefolds.
\endabstract
\endtopmatter
\define\codim{\operatorname{codim}}
\define\Supp{\operatorname{Supp}}
\define\Bs{\operatorname{Bs}}
\define\Exc{\operatorname{Exc}}
\define\Diff{\operatorname{Diff}}
\define\Rat{\operatorname{Rat}}
\document

\head
1. Introduction
\endhead

In this paper every variety is proper over the field $\Bbb C$ of complex numbers.
We follow the notation and terminology of \cite{Utah}.

Let $X$ be a normal algebraic variety and $\Delta = \sum d_i \Delta_i$ a $\Bbb 
Q$-divisor with $0 \le d_i \le 1$ on $X$ such that the log canonical divisor 
$K_X+\Delta$ is $\Bbb Q$-Cartier.
We call $(X,\Delta)$ a log pair.

Let $D$ be a nef (numerically effective) $\Bbb Q$-Cartier $\Bbb Q$-divisor on $X$.
We define the numerical Iitaka dimension $\nu (X,D):=\max\{ e;\ (D^e,S)>0$ for some 
subvariety $S$ of dimension $e$ on $X \}$.
The divisor $D$ is {\it abundant} if the Iitaka dimension $\kappa(X,D)$ equals $\nu (X,
D)$.
If, for some positive integer $m$, the divisor $mD$ is Cartier and the linear 
system $\vert mD \vert$ is free from base points, $D$ is said to be {\it 
semi-ample}.

For a birational morphism $f:Y \to X$ between normal algebraic varieties and for a 
divisor $E$ on $X$, the symbol $f^{-1}_*E$ expresses the strict transform of $E$ by 
$f$ and $f^{-1}(E)$ the set-theoretical inverse image.
A resolution $\mu: Y \to X$ is said to be a {\it log resolution} of the log pair $(X,\Delta)$ if the support of the divisor $\mu^{-1}_* \Delta + \sum \{ E; E \text{ is a } \mu \text{-exceptional prime}$ $\text{divisor} \}$ is with only simple normal crossings.
The log pair $(X,\Delta)$ is {\it log terminal} if there exists a log resolution 
$\mu :Y \to X$ such that $K_Y+\mu^{-1}_*\Delta=\mu^*(K_X+\Delta)+\sum a_i E_i$ with 
$a_i > -1$.
Moreover, if $\Exc (\mu)$ consists of divisors, $(X,\Delta)$ is said to be {\it 
divisorial log terminal} ({\it dlt}).
Szab\'o (\cite{Sz}) proved that the notions of dlt and wklt in \cite{Sh} are 
equivalent.
In the case where $(X,\Delta)$ is log terminal and $\lfloor \Delta \rfloor =0$, we 
say that $(X,\Delta)$ is {\it Kawamata log terminal} ({\it klt}).

We note that if $(X,\Delta)$ is klt then it is dlt.
In the Iitaka classification theory of open algebraic varieties, one embeds a 
smooth affine variety $U$ in some smooth projective variety $X$ such that $X 
\setminus U = \Supp (\Delta)$ where $\Delta$ is a reduced simple normal crossing 
divisor and studies the log pair $(X,\Delta)$.
In this case $(X,\Delta)$ is not klt but dlt.
Moreover it is known that we have to work allowing the $\Bbb Q$-factorial dlt singularities, to execute the log minimal model program for open algebraic varieties (see \cite{KMM}).
Therefore it is valuable to extend theorems proved in the case of klt pairs to the 
case of dlt pairs.

Now, concerning the log minimal model program, we review the famous

\proclaim{Log Abundance Conjecture (cf.\ \cite{KeMaMc})}
Assume that X is projective and $(X,\Delta)$ is dlt. If 
$K_X+\Delta$ is nef, then $K_X+\Delta$ is semi-ample.
\endproclaim

This conjecture claims that the concept of ``log minimal'' (that is, the log canonical divisor is nef) should be not only numerical but also geometric.
Kawamata (\cite{Ka1}) and Fujita (\cite{Fujt}) proved the conjecture in $\dim X =2$ 
and Keel, Matsuki and McKernan (\cite{KeMaMc}) in $\dim X =3$.
(The assumption concerning singularities in their papers is that $(X, \Delta)$ is log canonical, which is more general than dlt.)
Moreover Fujino proved

\proclaim{Theorem 1 (\cite{Fujn2, 3.1})}
Assume that $(X,\Delta)$ is dlt and $\dim X=4$.
If $K_X+\Delta$ is nef and big, then $K_X+\Delta$ is semi-ample.
\endproclaim

The following two theorems due to Kawamata are helpful to deal with the conjecture.

\proclaim{Theorem 2 (\cite{Ka2, 6.1})}
Assume that $(X,\Delta)$ is klt and $K_X+\Delta$ is nef.
If $K_X+\Delta$ is abundant, then it is semi-ample.
\endproclaim

\proclaim{Theorem 3 (\cite{Ka2, 7.3}, cf.\ \cite{KeMaMc, 5.6})}
Assume that $(X,\Delta)$ is klt and $K_X+\Delta$ is nef.
If $\kappa (X,K_X+\Delta) >0$ and the log minimal model and the log abundance 
conjectures hold in dimension $\dim X-\kappa (X,K_X+\Delta)$, then $K_X+\Delta$ is 
semi-ample.
\endproclaim

In this paper we try to generalize the above-mentioned theorems and obtain the 
following

\proclaim{Main Theorem}
Assume that $(X,\Delta)$ is dlt and $\dim X=4$.
If $K_X+\Delta$ is nef and $\kappa (X,K_X+\Delta) >0$, then $K_X+\Delta$ is 
semi-ample.
\endproclaim

We prove Main Theorem, along the lines in the proofs of Theorems 1 and 2, using 
Fujino's abundance theorem for semi log canonical threefolds which are not necessarily irreducible
(For the definition of the concept ``sdlt'' appearing below, see Definition 2 in Section 2.):

\proclaim{Theorem 4 (\cite{Fujn1})}
Let $(S,\Theta)$ be a sdlt threefold. If $K_S+\Theta$ is nef, then $K_S+\Theta$ is 
semi-ample.
\endproclaim

\remark{Remark}
If the log minimal model and the log abundance conjectures hold in dimension $\le 
n-1$, and Theorem 4 holds in dimension $n-1$, then Main Theorem holds in dimension 
$n$.
\endremark

{\it Acknowledgment.}
The Author would like to thank the referees for their valuable advice concerning the presentation and the quotations.

\head
2. Preliminaries
\endhead

In this section we state notions and results needed in the proof of Main Theorem.

The next two propositions are from the theories of the Kodaira-Iitaka dimension and 
the minimal model respectively.

\proclaim{Proposition 1 (\cite{Ii, Theorem 10.3})}
Let $D$ be an effective divisor on a smooth variety $Y$.
Suppose that the rational map $\Phi_{\vert D \vert}:Y \to Z$ is a morphism and that 
the rational function field $\Rat (\Phi_{\vert mD \vert}(Y))$ is isomorphic to 
$\Rat (Z)$ for all positive integer $m$.
Then $\Rat (Z)$ is algebraically closed in  $\Rat (Y)$ and $\kappa(W,D \vert_W)=0$ 
for a ``general'' fiber of $\Phi_{\vert D \vert}$ .
\endproclaim

\proclaim{Proposition 2 (\cite{KMM, Section 5-1.})}
Assume that $(X_{lm},\Delta_{lm})$ is a log minimal model for a $\Bbb Q$-factorial, 
dlt projective variety $(X,\Delta)$.
Then every common resolution $X\overset g\to\longleftarrow Y\overset 
h\to\longrightarrow X_{lm}$ satisfies the condition that $K_Y+g^{-1}_* \Delta +E 
\ge g^*(K_X+\Delta) \ge h^*(K_{X_{lm}}+\Delta_{lm})$, where $E$ is the reduced 
divisor composed of the $g$-exceptional prime divisors.
\endproclaim

The following is a vanishing theorem of Koll\'ar-type:

\proclaim{Theorem 5 (\cite{Ko, 10.13}, cf.\ \cite{Ka2, 3.2}, \cite{EV, 3.5})}
Let $f:X \to Y$ be a surjective morphism from a smooth projective variety $X$ to a 
normal variety $Y$.
Let $L$ be a divisor on $X$ and $D$ an effective divisor on $X$ such that $f(D) \ne 
Y$.
Assume that $(X,\Delta)$ is klt and $L-D-(K_X+\Delta)$ is $\Bbb Q$-linearly 
equivalent to $f^*M$ where $M$ is a nef and big $\Bbb Q$-Cartier $\Bbb Q$-divisor 
on $Y$.
Then the homomorphisms $H^i(X, \Cal O_X(L-D)) \to H^i(X, \Cal O_X(L))$
are injective for all $i$.
\endproclaim

When we work on the non-klt locus $\lfloor \Delta \rfloor$ of 
a log terminal pair $(X,\Delta)$, we need

\proclaim{Lemma 1 (cf.\ \cite{Ii, Proposition 1.43})}
Let $S$ be a reduced scheme and $\Cal F$ an invertible sheaf on $S$.
Then the restriction map $H^0(S,\Cal F) \to H^0(U,\Cal F)$ is injective for all 
open dense subset $U$ of $S$.
\endproclaim

The following lemma is used to manage cases where Theorem 5 can not be 
applied (See \cite{KeMaMc, Section 7}):

\proclaim{Lemma 2 (cf.\ \cite{Fujt0, 1.20})}
Let $f:S \to Z$ be a surjective morphism between normal varieties and $H_Z$ a 
Cartier divisor on $Z$.
If $f^*H_Z$ is semi-ample, then so is $H_Z$ .
\endproclaim

The set $\text{\bf Strata}(D)$ defined below is the set of non-klt centers for a smooth pair $(Y,D)$.

\definition{Definition 1}
Let $D=\sum_{i=1}^l D_i$ be a reduced simple normal crossing divisor on a smooth 
variety Y. We set $\text{\bf Strata}(D):=\{\Gamma ;\ 1 \le i_1<i_2<\dots <i_k \le l$, 
$\Gamma$ is an irreducible component of $ D_{i_1} \cap D_{i_2} \cap \dots \cap 
D_{i_k} \ne \emptyset \}$.
\enddefinition

When we manage the non-klt locus $\lfloor \Delta \rfloor$ of a dlt pair $(X,\Delta)
$, we need the following notion:

\definition{Definition 2 (due to Fujino (\cite{Fujn1, 1.1}))}
Let $S$ be a reduced $S_2$ scheme which is pure $n$-dimensional and normal crossing 
in dimension 1.
Let $\Theta$ be an effective $\Bbb Q$-Weil divisor such that $K_S+\Theta$ is $\Bbb 
Q$-Cartier.
Let $S=\bigcup S_i$ be the decomposition into irreducible components.
The pair $(S,\Theta)$ is {\it semi divisorial log terminal} ({\it sdlt}) if $S_i$ 
is normal and $(S_i,\Theta \vert_{S_i})$ is dlt for all $i$.
\enddefinition

\proclaim{Proposition 3 (\cite{Fujn1, 1.2.(3)}, cf.\ \cite{Sh, 3.2.3}, \cite{KoM, 
5.52})}
If $(X,\Delta)$ is dlt, then $(\lfloor \Delta \rfloor, \Diff (\Delta-\lfloor \Delta 
\rfloor))$ is sdlt.
\endproclaim

\head
3. Proof of Main Theorem
\endhead

The following proposition is used to imply the abundance of some log canonical divisor from its mobility:

\proclaim{Proposition 4 (\cite{Ka2, 7.3}, \cite{KeMaMc, 5.6})}
Let $(X,\Delta)$ be a variety with only log canonical singularities such that 
$K_X+\Delta$ is nef and $\kappa(X,K_X+\Delta) >0$.
If the log minimal model and the log abundance conjectures hold in dimension $\dim 
X-\kappa(X,K_X+\Delta)$, then $\kappa(X,K_X+\Delta) = \nu(X,K_X+\Delta)$.
\endproclaim

In the literature (Theorem 3 \cite{Ka2, 7.3}), this is proved for klt pairs.
However the proof is valid for log canonical pairs also.
Thus in the proof below we note only the parts where we have to be careful in reading \cite{Ka2, Proof of 7.3}.

\demo{Proof {\rm (}\cite{Ka2, Proof of 7.3 }{\rm )}}
By Proposition 1, we have a diagram $X\overset \mu\to\longleftarrow Y\overset 
f\to\longrightarrow Z$ with the following properties:
\roster
 \item"{(a)}" $Y$ and $Z$ are smooth projective varieties.
Moreover $Y$ is a log resolution of $(X,\Delta)$.
 \item"{(b)}" $\mu$ is birational and $f$ is surjective.
The morphism $f$ satisfies that $\dim Z=\kappa(X,K_X+\Delta)$ and $f_* \Cal 
O_Y=\Cal O_Z$ .
 \item"{(c)}" $K_Y+\mu^{-1}_* \Delta +E =\mu^*(K_X+\Delta)+E_{\mu}$ , where E is 
the reduced divisor composed of the $\mu$-exceptional prime divisors and $E_{\mu}$ 
is an effective $\Bbb Q$-divisor.
 \item"{(d)}" For a general fiber $W=Y_z$ of $f$, $K_Y \vert_W =K_W$ and $\kappa(W,
K_W+(\mu^{-1}_* \Delta +E)\vert_W)=0$.
\endroster
We note that $W$ is smooth and $\Supp ((\mu^{-1}_* \Delta +E)\vert_W)$ is with only 
simple normal crossings.

We apply the log minimal model program to $(W,(\mu^{-1}_* \Delta +E)\vert_W)$ and 
obtain a log minimal model $(W_{lm},\Delta_{lm})$, where $K_{W_{lm}}+\Delta_{lm} 
\sim_{\Bbb Q} 0$ from the log abundance.
We consider a common resolution $W\overset \rho\to\longleftarrow W'\overset 
\sigma\to\longrightarrow W_{lm}$ of $W$ and $W_{lm}$ such that $W'$ is projective.
From Proposition 2, 
$$\rho^*(K_W+(\mu^{-1}_* \Delta +E)\vert_W) =\sigma^* 
(K_{W_{lm}}+\Delta_{lm})+E_{\sigma} \sim_{\Bbb Q} E_{\sigma}$$
for some $\sigma$-exceptional effective $\Bbb Q$-divisor $E_{\sigma}$ .
Thus we have the relation
$$\rho^*(\mu^*(K_X+\Delta)\vert_W) = \rho^* (K_W + (\mu^{-1}_* \Delta +E)\vert_W 
-E_{\mu} \vert_W) \sim_{\Bbb Q} E_{\sigma} - {\rho}^*(E_{\mu} \vert_W).$$
We put $E_{+}-E_{-}:=E_{\sigma}-{\rho}^*(E_{\mu} \vert_W)$, where $E_{+}$ and 
$E_{-}$ are effective $\Bbb Q$-divisors that have no common irreducible components.
Here $E_{+}$ is $\sigma$-exceptional.

This paragraph is due to an argument in Miyaoka \cite{Mi, IV 2.4}.
Put $e:=\dim W'$ and $c:=$ the codimension of $\sigma(E_{+})$ in $W_{lm}$ .
We take general members $A_1,A_2,\dots,A_{e-c} \in \vert A \vert$ and 
$H_1,H_2,\dots,H_{c-2} \in \vert H \vert$ where $A$ and $H$ are very ample divisors 
on $W_{lm}$ and $W'$ respectively.
Set
$$S=(\bigcap_{i=1}^{e-c} \sigma^{-1}(A_i)) \cap (\bigcap_{i=1}^{c-2} H_i).$$

Taking into account the argument above, we proceed along the lines in \cite{Ka2, 
Proof of 7.3}.
Then we have the fact that $\rho^*(\mu^*(K_X+\Delta)\vert_W)$ is $\Bbb Q$-linearly 
trivial and so is $\mu^*(K_X+\Delta)\vert_W$ .
From this the assertion follows.
\qed
\enddemo

In the following we cope with the base points that lie on the non-klt locus $\lfloor \Delta \rfloor$:

\proclaim{Proposition 5}
Let $(X,\Delta)$ be a log terminal variety and H a nef $\Bbb Q$-Cartier $\Bbb 
Q$-divisor such that $H-(K_X+\Delta)$ is nef and abundant.
Assume that $\nu (X,aH-(K_X+\Delta)) = \nu (X,H-(K_X+\Delta))$ and $\kappa 
(X,aH-(K_X+\Delta)) \ge 0$ for some $a \in \Bbb Q$ with $a>1$.
If $H \vert _{\lfloor \Delta \rfloor}$ is semi-ample, then $\Bs \vert mH \vert \cap 
\lfloor \Delta \rfloor = \emptyset$ for some positive integer $m$ with $mH$ being 
Cartier.
\endproclaim

\demo{Proof}
From an argument in \cite{Ka2, Proof of 6.1} we have a diagram $X\overset \mu\to\longleftarrow Y\overset f\to\longrightarrow Z$ 
with the following properties:
\roster
 \item"{(a)}" $Y$ and $Z$ are smooth projective varieties.
Moreover $Y$ is a log resolution of $(X,\Delta)$.
 \item"{(b)}" $\mu$ is birational and $f$ is surjective with the property that $f_* 
\Cal O_Y = \Cal O_Z$ .
 \item"{(c)}" $\mu^*(H-(K_X+\Delta)) \sim_{\Bbb Q} f^*M_0$ for some nef and big 
$\Bbb Q$-divisor $M_0$ .
 \item"{(d)}" $\mu^*H \sim_{\Bbb Q} f^*H_0$ for some nef $\Bbb Q$-divisor $H_0$ .
\endroster
We define rational numbers $a_i$ by $K_Y=\mu^*(K_X+\Delta)+\sum a_i E_i$ .
We may assume that $H_0$ and $H$ are Cartier.

We put
$$S:=\lfloor \Delta \rfloor,\ E:=\sum_{a_i>0}\lceil a_i \rceil E_i \ \text{and} \ 
S':=\sum_{a_i=-1} E_i \ .$$
We note that $m\mu^*H+E-S'-(K_Y+\sum\{ -a_i \} 
E_i)=(m-1)\mu^*H+\mu^*(H-(K_X+\Delta))$, which is $\Bbb Q$-linearly equivalent to 
the inverse image of a nef and big $\Bbb Q$-divisor on $Z$.
There are two cases:
\enddemo

\demo{Case {\rm 1}}
$f(S') \ne Z$.
\quad
In this case we use Fujino's argument \cite{Fujn2, Section 2}.
By Theorem 5 we have an injection
$$H^1(Y,\Cal O_Y(m\mu^*H+E-S')) \to H^1(Y,\Cal O_Y(m\mu^*H+E)).$$
Then we consider the commutative  diagram:
$$\CD
 H^0(Y,\Cal O_Y(m\mu^*H+E)) @>\text{surjective}>> H^0(S',\Cal O_{S'}(m\mu^*H+E)) 
@>>> 0 \\
 @A{\cong}AA @AAiA \\
 H^0(Y,\Cal O_Y(m\mu^*H)) @>>> H^0(S',\Cal O_{S'}(m\mu^*H))\\
 @A{\cong}AA @AAjA @. \\
 H^0(X,\Cal O_X(mH)) @>s>> H^0(S,\Cal O_S(mH)) @.
\endCD$$
The homomorphism $i$ is injective from the fact that $E$ and $S'$ have no common 
irreducible component and Lemma 1.
The homomorphism $j$ is injective from the fact that $S' \to S$ is surjective and 
Lemma 1.
Thus the homomorphism $s$ is surjective from the diagram.
Consequently $\vert mH \vert \big\vert_S = \big\vert mH \vert_S \big\vert$.
\enddemo

\demo{Case {\rm 2}}
$f(S')=Z$.
\quad
In this case we use an argument in \cite{KeMaMc, Section 7}.
There exists an irreducible component $S''$ of $S'$ such that $f(S'')=Z$.
Because $H \vert_S$ is semi-ample and $\mu^*H \sim_{\Bbb Q} f^*H_0$ , $f^*H_0 
\vert_{S''}$ is semi-ample.
Consequently the $\Bbb Q$-divisor $H_0$ also is semi-ample from Lemma 2.
\qed
\enddemo

We generalize Kawamata's result \cite{Ka2, 6.1} (see also Theorem 2) concerning the semi-ampleness for klt pairs to the case of log terminal pairs in the following form:

\proclaim{Proposition 6}
Assume that $(X,\Delta)$ is log terminal.
Let $H$ be a nef $\Bbb Q$-Cartier $\Bbb Q$-divisor on $X$ with the following 
properties:
\roster
 \item $H-(K_X+\Delta)$ is nef and abundant.
 \item $\nu (X,aH-(K_X+\Delta)) = \nu (X,H-(K_X+\Delta))$ and $\kappa 
(X,aH-(K_X+\Delta)) \ge 0$ for some $a \in \Bbb Q$ with $a>1$.
\endroster
If, for some positive integer $p_1$ , the divisor $p_1H$ is Cartier and $\Bs \vert 
p_1H \vert \cap \lfloor \Delta \rfloor = \emptyset$, then $H$ is semi-ample.
\endproclaim

In the proof below we proceed along the lines in \cite{Ka2, Proof of 6.1} and thus omit the parts which are parallel.
However we have to be very delicate in dealing with the non-klt locus $\lfloor \Delta \rfloor$.

\demo{Proof}
From \cite{Ka2, 6.1}, we may assume that $\lfloor \Delta \rfloor \ne 0$.
Therefore the condition that $\Bs \vert p_1 H \vert \cap \lfloor \Delta \rfloor = \emptyset$ implies that $\Bs \vert p_1 H \vert \ne X$.
Thus $\vert p_1tH \vert \ne \emptyset$ for all $t \in \Bbb N_{>0}$ (where $\Bbb 
N_{>0}$ denotes the set of all positive integers).

We have smooth projective varieties $Y$ and $Z$ and morphisms $X\overset 
\mu\to\longleftarrow Y\overset f\to\longrightarrow Z$ with the following 
properties:
\roster
 \item $\mu$ is birational and $f$ is surjective.
 \item $f_* \Cal O_Y = \Cal O_Z$ .
 \item $\mu^*(H-(K_X+\Delta)) \sim_{\Bbb Q} f^*M_0$ for some nef and big $\Bbb 
Q$-divisor $M_0$ (where the symbol $\sim_{\Bbb Q}$ expresses the $\Bbb Q$-linear 
equivalence).
 \item $\mu^*H \sim_{\Bbb Q} f^*H_0$ for some nef $\Bbb Q$-divisor $H_0$ .
\endroster
We may assume that $H_0$ and $H$ are Cartier and $f^*H_0$ and $\mu^*H$ are linearly 
equivalent.

Putting $\Lambda (m):= \Bs \vert mH \vert$, we may assume that $\Lambda (p_1) \ne 
\emptyset$ (otherwise we immediately obtain the assertion).
By repetition of blowing-ups over $Y$, we may replace $Y$ and get a simple normal 
crossing divisor $F=\sum_{i \in I}F_i$ on $Y$ such that 
\roster
 \item"{(5)}" $\mu^*\vert p_1H \vert = \vert L \vert +\sum_{i \in I} r_iF_i$ and 
$\vert L \vert$ is base point free.
\endroster
Then by replacing $Z$ and $Y$ we have $L \sim_{\Bbb Q} f^*L_0$ for some $\Bbb 
Q$-divisor $L_0$ , because $\nu(Y,\mu^*(aH-(K_X+\Delta))) \ge \nu(Y,((a-1)/p_1)L+  
\mu^*(H-(K_X+\Delta))) \ge \nu(Y,\mu^*(H-(K_X+\Delta)))$ from the argument in 
\cite{Ka2, Proof of 2.1}.
We note that
$$\Lambda (p_1) = \mu (\bigcup_{r_i \ne 0} F_i).$$

We have an effective divisor $M_1$ such that $M_0-\delta M_1$ is ample for all 
$\delta \in \Bbb Q$ with $0<\delta \ll 1$.
By further repetition of blowing-ups over $Y$, we may replace $Y$ and get the 
following properties:
\roster
 \item"{(6)}" $K_Y = \mu^*(K_X+\Delta)+\sum_{i \in I} a_iF_i$ .
 \item"{(7)}" $f^*M_1=\sum_{i \in I} b_iF_i$ .
\endroster

We set
$$c:=\min_{r_i \ne 0} \frac{a_i+1-\delta b_i}{r_i}.$$
Note that if $a_i=-1$ then $\mu (F_i) \subset \lfloor \Delta \rfloor$ and that if 
$\mu (F_i) \subset \lfloor \Delta \rfloor$ then $r_i=0$ from the assumption of the 
theorem.
Thus by taking $\delta$ small enough, we may assume that $c>0$ and that, if $F_i 
\not\subset \mu^{-1}(\lfloor \Delta \rfloor$), then $a_i+1-\delta b_i >0$ (even if 
$b_i \ne 0$).
Set $I_0:=\{i \in I;\ a_i+1-\delta b_i=cr_i,\ r_i \ne 0 \}$ and $\{ Z_{\alpha} 
\}:=\{f(\Gamma);\ \Gamma \in \text{\bf Strata} (\sum_{i \in I_0}F_i)\}$.
Let $Z_1$ be a minimal element of $\{Z_{\alpha}\}$ with respect to the inclusion 
relation.
We note that $Z_1 \ne Z$.
Because $M_0-\delta M_1$ is ample, for some $q \in \Bbb N_{>0}$ , there exists a 
member $M_2 \in \vert q(M_0-\delta M_1) \vert$ such that $Z_1 \subset M_2$ and 
$Z_{\alpha} \not\subset M_2$ for all $\alpha \ne 1$.

We would like to show that we may assume that $\Supp (f^* M_2) \subset F$.
Then we investigate the variation of the numbers $a_i+1-\delta b_i$ and the set 
$I_0$ under the blowing-up $\sigma :Y^{\prime} \to Y$ with permissible smooth 
center $C$ with respect to F.
We get a simple normal crossing divisor $F^{\prime} =\sum_{i \in I^{\prime}} 
F^{\prime}_i$ on $Y^{\prime}$ (where $I'=I \cup \{ 0 \}$) with the following 
properties:
$$\align
 F^{\prime}_0 & = \sigma^{-1}(C).\tag 8\\
 K_{Y^{\prime}} & =\sigma^* \mu^* (K_X+\Delta)+\sum_{i \in I'} a^{\prime}_i 
F^{\prime}_i \ .\tag 9\\
 \sigma^*(\sum_{i \in I} r_i F_i) & = \sum_{i \in I'} r^{\prime}_i F^{\prime}_i \ .
\tag 10\\
 \sigma^* f^* M_1 & = \sum_{i \in I'} b^{\prime}_i F^{\prime}_i \ .\tag 11
\endalign$$
We set $I_0':=\{i \in I';\ a'_i+1-\delta b'_i=cr'_i,\ r'_i \ne 0 \}$.
Let $F_{i_1}, \dots, F_{i_u}$ be the irreducible components of $F$ that contain 
$C$.
Let $F'_{i_j}$ be the strict transform of $F_{i_j}$ by $\sigma$.
We note that
$$\sigma^*(K_Y-\sum_{j=1}^u a_{i_j} F_{i_j})=K_{Y'}-(\codim_Y C-1)F'_0 
-\sum_{j=1}^u a_{i_j} (F'_{i_j}+F'_0).$$
Thus $a'_0=(\codim_Y C-1)+\sum_{j=1}^u a_{i_j}$ .
Therefore
$$a'_0+1 \ge \sum_{j=1}^u (a_{i_j}+1), \tag 12$$
where the equality holds if and only if $u=\codim_Y C$.
We note also that $r'_0=\sum_{j=1}^u r_{i_j}$ and $b'_0=\sum_{j=1}^u b_{i_j}$ .
\enddemo

\proclaim{Claim 1}
If $F'_0 \not\subset (\mu \sigma)^{-1}(\lfloor \Delta \rfloor)$, then 
$a'_0+1-\delta b'_0 \ge cr'_0$ .
The equality holds if and only if $\codim_Y C =u$ and $i_j \in I_0$ for all $j$.
\endproclaim

\demo{Proof of Claim {\rm 1}}
First  we note the inequality
$$a'_0+1-\delta b'_0 \ge \sum _{j=1}^u (a_{i_j}+1-\delta b_{i_j}),$$
where the equality holds if and only if $\codim_Y C =u$.
Because $F_{i_j} \not\subset \mu^{-1}(\lfloor \Delta \rfloor)$, we have 
$a_{i_j}+1-\delta b_{i_j} >0$.
Here if $r_{i_j} \ne 0$ then $a_{i_j}+1-\delta b_{i_j} \ge c r_{i_j}$ , from the 
definition of $c$.
On the other hand if $r_{i_j} = 0$ then $a_{i_j}+1-\delta b_{i_j} > c r_{i_j}$ .
Now we note the inequality
$$\sum _{j=1}^u (a_{i_j}+1-\delta b_{i_j}) \ge \sum _{j=1}^u cr_{i_j}\ ,$$
where the equality holds if and only if $r_{i_j} \ne 0$ and $a_{i_j}+1-\delta 
b_{i_j} = c r_{i_j}$ (that is, $i_j \in I_0$) for all $j$.
Here $\sum _{j=1}^u cr_{i_j} =cr'_0$ .
{\it Proof of Claim 1 ends.}
\enddemo

\proclaim{Claim 2}
If $i_j \in I_0$ for all $j$ and $C \in \text{\bf Strata} (\sum_{j=1}^u F_{i_j})$, 
then $I_0' = I_0 \cup \{ 0 \}$.
Otherwise $I_0' = I_0$ .
\endproclaim

\demo{Proof of Claim {\rm 2}}
We note that $\codim_Y C = u$ if and only if $C \in \text{\bf Strata} (\sum_{j=1}^u 
F_{i_j})$.
Thus Claim 1 implies the assertion, because if $F'_0 \subset (\mu 
\sigma)^{-1}(\lfloor \Delta \rfloor)$ then $r'_0 =0$.
{\it Proof of Claim 2 ends.}
\enddemo

\proclaim{Claim 3}
$$\min_{r'_i \ne 0} \frac{a'_i+1-\delta b'_i}{r'_i} = c.$$
\endproclaim

\demo{Proof of Claim {\rm 3}}
In the case where $r'_0 \ne 0$, we have $F'_0 \not\subset (\mu \sigma)^{-1}(\lfloor 
\Delta \rfloor)$.
Thus Claim 1 implies the assertion.
{\it Proof of Claim 3 ends.}
\enddemo

\proclaim{Claim 4}
If $F'_0 \not\subset (\mu \sigma)^{-1}(\lfloor \Delta \rfloor)$, then 
$a'_0+1-\delta b'_0 >0$.
\endproclaim

\demo{Proof of Claim {\rm 4}}
In this case, $a'_0+1 >0$.
If $b'_0 \ne 0$, then $C \subset f^*M_1$ , so $u \ne 0$.
Thus $a'_0+1-\delta b'_0 \ge \sum_{j=1}^u (a_{i_j}+1-\delta b_{i_j}) >0$ because 
all $F_{i_j} \not\subset \mu^{-1}(\lfloor \Delta \rfloor)$.
{\it Proof of Claim 4 ends.}
\enddemo

{\it Proof of Proposition 6 continues.}
By virtue of Claims 2, 3 and 4, we may assume that $f^*M_2 = \sum_{i \in I} s_i 
F_i$ where $F = \sum_{i \in I} F_i$ is a simple normal crossing divisor.
We put
$$c':=\min_{\mu(F_i) \not\subset \lfloor \Delta \rfloor} 
\frac{a_i+1-\delta b_i}{r_i+\delta' s_i} \tag 13$$
and $I_1:=\{ i \in I;\ a_i+1-\delta b_i=c'(r_i+\delta' s_i),\ \mu(F_i) \not\subset 
\lfloor \Delta \rfloor \}$, for a rational number $\delta'$ with $0<\delta' \ll 
\delta$.

\proclaim{Claim 5}
$I_1 \subset I_0$ .
\endproclaim

\demo{Proof of Claim {\rm 5}}
Because if $\mu(F_i) \not\subset \lfloor \Delta \rfloor$ then $a_i+1-\delta b_i 
>0$, in the case where $r_i = 0$ the divisor $F_i$ does not attain the minimum in 
(13).
{\it Proof of Claim 5 ends.}
\enddemo

\proclaim{Claim 6}
There exists a member $j \in I_0$ such that $s_j >0$.
\endproclaim

\demo{Proof of Claim {\rm 6}}
The condition that $Z_1 \subset M_2$ implies that, for some $j \in I$, $s_j>0$ and 
$F_j$ contains an element $\Gamma \in \text{\bf Strata} (\sum_{i \in I_0}F_i)$.
Here $j \in I_0$ , because $F$ is with only simple normal crossings.
{\it Proof of Claim 6 ends.}
\enddemo

\proclaim{Claim 7}
$s_i >0$ for all $i \in I_1$ .
\endproclaim

\demo{Proof of Claim {\rm 7}}
Claims 5 and 6 and the formula (13) imply the assertion.
{\it Proof of Claim 7 ends.}
\enddemo

\proclaim{Claim 8}
$f(\Gamma)=Z_1$ for all $\Gamma \in \text{\bf Strata} (\sum_{i \in I_1}F_i)$.
\endproclaim

\demo{Proof of Claim {\rm 8}}
From Claim 7, $f(\Gamma) \subset M_2$ .
The condition that $Z_{\alpha} \not\subset M_2$ for all $\alpha \ne 1$ implies the 
fact that $f(\Gamma) \ne Z_{\alpha}$ for all $\alpha \ne 1$.
Thus $f(\Gamma)= Z_1$ from Claim 5.
{\it Proof of Claim 8 ends.}
\enddemo

{\it Proof of Proposition 6 continues.}
Now we set $N:=m \mu^* H + \sum_{i \in I} (-c'(r_i+\delta' s_i)+a_i-\delta 
b_i)F_i-K_Y$ for an integer $m \ge c'p_1+1$.
Then
$$\align
 N & =c'(-\sum_{i \in I} r_i F_i)+m \mu^*H -\mu^*(K_X+\Delta)-\delta\sum_{i \in I} 
b_i F_i -c'\delta'\sum_{i \in I} s_i F_i \\
   & \sim_{\Bbb Q} c'(L-p_1 \mu^*H)+m \mu^*H -\mu^*H +f^*(M_0-\delta 
M_1)-c'\delta'\sum_{i \in I} s_i F_i \\
   & \sim_{\Bbb Q} c'f^*L_0+(m-(c'p_1+1))\mu^*H +(1-c'\delta'q)f^*(M_0-\delta M_1).
\endalign$$
Because $\mu^*H$ and $f^*H_0$ are linearly equivalent, $N$ is $\Bbb Q$-linearly 
equivalent to the pull back of an ample $\Bbb Q$-divisor on $Z$.
We put
$$\align
 A: & =\sum_{i \in I \setminus I_1 \ \text{and} \ \mu (F_i) \not\subset \lfloor 
\Delta \rfloor} (-c'(r_i+\delta' s_i)+a_i-\delta b_i)F_i \ ,\\
 B_1: & =\sum_{i \in I_1} F_i \ ,\\
 C: & =\sum_{\mu (F_i) \subset \lfloor \Delta \rfloor} (-c'(r_i+\delta' 
s_i)+a_i-\delta b_i)F_i \ .
\endalign$$
Then $\sum_{i \in I} (-c'(r_i+\delta' s_i)+a_i-\delta b_i)F_i =A-B_1+C$.
We express $\lceil C \rceil:=-B_2+B_3$ in effective divisors $B_2$ and $B_3$ 
without common irreducible components.
Here we note that $f(B_1+B_2) \ne Z$, from Claim 8 and from the fact that the locus 
$f^{-1} (\Bs \vert p_1 H_0 \vert) = \mu^{-1} (\Lambda(p_1)) \ne \emptyset$ and the 
locus $\mu^{-1} (\lfloor \Delta \rfloor)$ are mutually disjoint.
Note also that $\lceil A \rceil$ and $B_3$ are $\mu$-exceptional effective divisors 
because if $a_i >0$ then $F_i$ is $\mu$-exceptional.

By Theorem 5, the homomorphism
$$H^1(Y,\Cal O_Y(m \mu^*H+\sum_{i \in I} \lceil (-c'(r_i+\delta' s_i)+a_i-\delta 
b_i) \rceil F_i)) \to H^1(Y,\Cal O_Y(m \mu^*H+ \lceil A \rceil + B_3))$$ is 
injective because $f(B_1+B_2) \ne Z$.
Hence 
$$\align
 & H^0(Y,\Cal O_Y(m \mu^*H+ \lceil A \rceil + B_3))\\
 & \quad \to H^0(B_1,\Cal O_{B_1}(m \mu^*H+ \lceil A \rceil + B_3)) \oplus H^0(B_2,
\Cal O_{B_2}(m \mu^*H+ \lceil A \rceil + B_3))
\endalign$$
is surjective, because $B_1 \cap B_2 = \emptyset$ from Claim 5.
Here
$$H^0(B_1,\Cal O_{B_1}(m \mu^*H+ \lceil A \rceil + B_3)) \cong H^0(B_1,\Cal O_{B_1}
(m \mu^*H+ \lceil A \rceil))$$
because $B_1 \cap B_3 = \emptyset$ from Claim 5.
We note that $\Supp(A \vert _{B_1})$ is with only simple normal crossings and 
$\lceil A \vert _{B_1} \rceil$ is effective.
Because $m\mu^*H \vert_{B_1} + A \vert _{B_1}-K_{B_1} = N \vert _{B_1}$ , we obtain 
a positive integer $p_2$ such that
$$H^0(B_1,\Cal O_{B_1}(p_2 t \mu^*H+ \lceil A \rceil)) \ne 0$$
for all $t \gg 0, $ from Claim 8 and \cite{Ka2, 5.1}.
Consequently the assertion of the proposition follows.
\qed

\demo{{\bf Proof of Main Theorem}}
Because $\kappa(X,K_X+\Delta)>0$, we have $\kappa(X,K_X+\Delta) =\nu(X,K_X+\Delta)$ 
from the log minimal model and the log abundance theorems in dimension $\le3$ 
(\cite{Sh}, \cite{KeMaMc}) and Proposition 4.
We note that $(K_X+\Delta)\vert_{\lfloor \Delta \rfloor}$ is semi-ample from 
Proposition 3 and Theorem 4.
Thus Proposition 5 implies that $\Bs \vert m(K_X+\Delta) \vert \cap \lfloor \Delta 
\rfloor = \emptyset$ for some $m \in \Bbb N_{>0}$ with $m(K_X+\Delta)$ being 
Cartier.
Consequently Proposition 6 gives the assertion.
\qed
\enddemo

\enddocument